\documentclass[12pt]{article}
\usepackage{graphicx}
\usepackage{amsfonts}
\usepackage{bbm}
\usepackage{mathrsfs}
 \usepackage{mathrsfs,amsmath,amssymb,amsthm}
\usepackage{color}
 \usepackage{amssymb}
 \usepackage{amsbsy}
\allowdisplaybreaks
 \theoremstyle{plain}
\theoremstyle{remark}  \newtheorem{remark}{\noindent\mbox{Remark}}
 \theoremstyle{plain}
 \theoremstyle{plain}\newtheorem{lemma}{\noindent\mbox{Lemma}}
\theoremstyle{plain} \newtheorem{theorem}{\noindent\mbox{Theorem}}
 \theoremstyle{plain}
 \theoremstyle{plain}
\theoremstyle{definition}

 \def\bq{\begin{equation}}
 \def\eq{\end{equation}}
 \def\eqn{\end{eqnarray}}
 \def\bqn{\begin{eqnarray}}
 \def\proof{\noindent{\it Proof.~~}}
 \def\qed{\hfill$\Box$\medskip}
 \def\rto{\rightarrow\infty}
 \def\z{\left}
 \def\y{\right}
 \def\no{\nonumber}
\def\mbz{\mathbb{Z}}
 \def\mbe{\mathbb{E}}
 \def\mbp{\mathbb{P}}

\begin{document}
\begin{center}{\Large \bf Range of $(1,2)$ random walk in random environment}\footnote{Supported by National
Nature Science Foundation of China (Grant No. 11501008),  Nature Science Foundation of Anhui Province (Grant No. 1508085QA12) and Nature Science Foundation of Anhui Educational Committee (Grant No. KJ2014A085).}\end{center}

\vspace{0.5cm}
\centerline {Hua-Ming \uppercase{Wang}$^a$ }

    \begin{center}
      {\footnotesize$^a$Department of Statistics, Anhui Normal University, Wuhu 241003, P. R. China

    E-mail\,$:$ hmking@mail.ahnu.edu.cn}
    \end{center}

\vspace{-.3cm}

\begin{center}
\begin{minipage}[c]{12cm}
\begin{center}\textbf{Abstract}\quad \end{center}

Consider $(1,2)$ random walk in random environment $\{X_n\}_{n\ge0}.$ In each step, the walk jumps at most a distance $2$ to the right or a distance $1$ to the left. For the walk transient to the right, it is proved that almost surely $\lim_{x\rightarrow\infty}\frac{\#\{X_n:\ 0\le X_n\le x,\ n\ge0\}}{x}=\theta$ for some $0<\theta<1.$ The result shows that the range of the walk covers only a linear proportion of the lattice of the positive half line. For the nearest neighbor random walk in random or non-random environment, this phenomenon could not appear in any circumstance.

\vspace{0.2cm}


\vspace{0.2cm}
\end{minipage}
\end{center}

\section{Introduction}
\subsection{Motivation}
It is well known that Random Walks in Random Environments (RWRE) own a lot of surprising phenomenons compared with simple random walk. In this paper we study one dimensional RWRE with bounded jumps which will be written as $(L,R)$ RWRE. In each step, the walk jumps at most a distance $R$ to the right and at most a distance $L$ to the left.
Basically speaking, the limiting behaviors of RWRE with bounded jumps are more or less the same as the nearest neighbor RWRE.
For example, both of them exhibit a slowdown property. Precisely, letting $\{X_n\}_{n\ge0}$ be the RWRE, then it is possible that $\lim_{n\rto}X_n=\infty$ but $\lim_{n\rto}\frac{X_n}{n}=0.$ For the details, see Solomon \cite{sol} for the nearest neighbor setting and  Key \cite{key} for the bounded-jumping setting.
Also they shares a so-called $\log^2$ law, that is, for the recurrent case, under certain conditions, $\frac{X_n}{\log^2n}\overset{d}{\rightarrow}X, $ for some non degenerate random variable $X,$   as $n\rto.$ See Sinai \cite{sinai} and Letchikov \cite{letc} for details.

We want to reveal some distinct limiting behavior caused by the so-called bounded jumps. By intuition, since the walk is non nearest neighbor, for the transient case, it may skip some points and go directly to the infinity. In this paper, letting $L=1, R=2,$ we study $(1,2)$ RWRE to illustrate this new phenomenon. We prove that for the $(1,2)$ RWRE which is transient to the right, the range of the walk covers just a linear proportion of $\mathbb Z_+:=\{0,1,2,...\}.$

 In literatures, for random walk in non random environment, $\#\{X_k:0\le k\le n\},$ the range up to time $n,$  was usually considered. The limit $$\lim_{n\rto}\frac{\#\{X_k:0\le k\le n\}}{n}$$  was always of the concern, where and throughout, we use notation $``\#\{\ \}"$ to count the number of elements in a set $``\{\ \}".$

In this paper, from a different point of view, we consider the range of $(1,2)$ RWRE $\{X_n\}_{n\ge0}.$
For $x>0,$ consider $\#\{X_n: 0\le X_n\le x, n\ge 0\},$  which counts the number of all sites in $[0, x]$ the walk has ever visited. We show that when the walk is transient to the right, \begin{equation*}\label{lp}\lim_{x\rto}\frac{\#\{X_n: 0\le X_n\le x, n\ge 0\}}{x}=\theta \text{ for some } \theta\in (0,1).\end{equation*}
This results shows that a linear proportion of $\mathbb Z_+$ is not visited by the walk.

Next we introduce precisely the model and state the main results.
\subsection{Model and result}
For $i\in \mathbb Z,$ let $\omega_i=(\omega_i(-1),\omega_i(1),\omega_i(2))$ be a probability measure on $\{i-1,i+1,i+2\}.$
Let $\Omega$ be the collection of all $\omega=(\omega_i)_{i\in\mathbb Z}.$ Equip $\Omega$ with the weak topology and let $\mathcal F$ be the Borel $\sigma$-algebra. Let $\mbp$ be a probability measure on $(\Omega,\mathcal F)$ which makes $\omega=(\omega_{i})_{i\in \mbz}$ an i.i.d. sequence. For a typical realization of $\omega,$ we consider a Markov chain $\{X_n\}_{n\ge 0}$ with transitional probabilities
\begin{equation*}\left\{\begin{array}{l}
                         P_{\omega}^{x_0}(X_{n+1}=i-1\big |X_{n}=i)=\omega_{i}(-1),\\
                         P_{\omega}^{x_0}(X_{n+1}=i+1\big |X_{n}=i)=\omega_{i}(1),\\
                          P_{\omega}^{x_0}(X_{n+1}=i+2\big |X_{n}=i)=\omega_{i}(2),\\
                       P_\omega^{x_0}(X_0=x_0)=1,
                      \end{array}
\right.\end{equation*}
so that $P_\omega^{x_0}$ is the quenched law of the Markov chain starting from $x_0$ in the
environment $\omega.$
Define a new probability measure $P^{x_0}$ by $$P^{x_0}(\cdot)=\int P^{x_0}_\omega(\cdot)\mbp(d\omega),$$ which is called the annealed probability. We use $E_\omega^{x_0},$ $E^{x_0}$ and $\mbe$ to denote the expectation operator for $P_\omega^{x_0},$ $P^{x_0}$ and $\mbp$ respectively. The superscript $x_0$ will be omitted whenever it is $0.$

For $i\in\mbz,$ set $a_i(1)=\frac{\omega_i(1)+\omega_{i}(2)}{\omega_{i}(-1)},$  $a_i(2)=\frac{\omega_{i}(2)}{\omega_{i}(-1)}$
and let $$A_i=\left(
               \begin{array}{cc}
                 a_i(1)& a_{i}(2) \\
                 1 & 0 \\
               \end{array}
             \right)
$$
Then $\{A_i\}_{i\in\mbz}$ is a sequence of i.i.d. random matrices under $\mbp.$

Note that under condition
\begin{equation}\label{ell}\mbe(\log\omega_0(-1))>-\infty,\ \mbe(\log\omega_0(2))>-\infty,\end{equation}
 $\mathbb E |\ln \|A_0^{-1}\||+\mathbb E|\ln \|A_0\||<\infty.$ Hence an application of Oseledec's multiplicative ergodic theorem (see \cite{osel}) to the sequence $\{A_i\}_{i\in \mathbb Z}$ yields the Lyapunov
exponents of the sequence $\{A_i\}_{i\in\mathbb Z}$ which we write in increasing order as
$$-\infty<\gamma_1\le \gamma_2<\infty.$$
  Furthermore  due to the positivity, the top Lyapunov exponent $$\gamma_2=\lim_{n\rto}\frac{1}{n}\log\|A_0A_1\cdots A_{n-1}\|,\ \mbp\text{-a.s.}.$$

  The following result could be found in Letchikov \cite{letca}.

 {\noindent \bf Recurrence criterion:}  Assume condition (\ref{ell}). Then

a) $\gamma_2>0\Rightarrow P(\lim_{n\rto}X_n=\infty)=1;$

b) $\gamma_2=0\Rightarrow P(-\infty=\liminf_{n\rto}X_n<\limsup_{n\rto}X_n=\infty)=1;$

c) $\gamma_2<0\Rightarrow P(\lim_{n\rto}X_n=-\infty)=1.$

We remark that in \cite{letca} the  Lyapunov exponents of $\{\overline A_i\}_{i\in\mbz}$ are used to give the recurrence criteria, where $$\overline A_i=\left(
               \begin{array}{cc}
                 0 & 1\\
                 a_{i}(2) & a_i(1) \\
               \end{array}
             \right).$$
This causes no problem because $\{A_i\}_{i\in \mbz}$ and $\{\overline A_i\}_{i\in \mbz}$ share the same Lyapunov exponents.
We are now ready to state the main theorem.
\begin{theorem}\label{main}
Suppose that condition (\ref{ell}) holds and $\gamma_2>0.$ Then $P$-a.s.,
  \begin{equation}\label{lpm}\lim_{x\rto}\frac{\#\{X_n: 0\le X_n\le x, n\ge 0\}}{x}=\theta \text{ for some } \theta\in (0,1).\no\end{equation}
\end{theorem}
A formula for $\theta$ will be given, see (\ref{taa}) (\ref{ta}) and (\ref{tb})  below. It is closely related to the maximum of a positive excursion of the walk.
Define
 $D=\inf\{n\ge 0:X_n< X_0\}$ and
 whenever $D<\infty$ set $$M=\sup\{X_n: 0\le n\le D\}$$ the maximum of a positive excursion,
 where and throughout  we use the conventions $\inf\emptyset=\infty$ and $\sup\emptyset=0.$
  We have the following results.
  \begin{theorem}\label{mt}
   Suppose that condition (\ref{ell}) holds and $\gamma_2>0.$ Then there exist $c_1,c_2>0$ such that
   $P(M>n,D<\infty)<c_1e^{-c_2n}$  for  all $n>0.$
 \end{theorem}

 To conclude the introduction, we give the idea of the proof. Theorem \ref{main} is proved by using a renewal structure. Let $0=\iota_0<\iota_1<\iota_2<...$ be the successive epochs of the walk. At each epoch $\iota_n$, the walk reaches a new level which exceeding the maximum of the walk before this epoch of a distance 2 and will never come back to any site below this new level. Therefore, the site between the new level and the former maximum is never visited by the walk. By showing that $\{X_{\iota_i}-X_{\iota_{i-1}}\}_{i\ge2}$ are independent random variables  distributed as $(X_{\iota_1})$ under $P(\cdot |D=\infty)$ and $E(X_{\iota_1}|D=\infty)<\infty,$ we conclude that a linear proportion of $\mbz_+$ is not visited by the walk. Theorem \ref{mt} is proved by a large deviation argument of the products of a sequence of i.i.d. random matrices to estimate the escaping probability of the walk.

\section{Renewal structure of $(1,2)$ RWRE}
In this section, we introduce a renewal structure which is crucial for us to prove Theorem \ref{main}. Up to our knowledge, this approach was first used by Kesten \cite{kes77} to prove a renewal theorem for the nearest neighbor RWRE, which was generalized to $(L,1)$ RWRE by Hong and Sun \cite{hs13}. We mention also that this approach was successfully used to study the high dimensional RWRE. For the details, see Sznitman \cite{sz04}, Sznitman and Zerner \cite{sz99} and reference therein. We borrow a lot of skills used in \cite{sz99} in the proof.

Define $D=\inf\{n\ge 0: X_n<X_0\}.$ Denote by $(\theta_n)_{n\ge 0}$ the canonical shift on $\mbz^{\mathbb N}.$  Let $S_0=0$ and $M_0=X_0.$ Define
\begin{align*}
  &S_1=\inf\{n:X_n>M_0\},\\
   &R_1=D\circ \theta_{S_1}+S_1,\\
   & M_1=\sup\{X_n:0\le n\le R_1\}
\end{align*}
and by induction, for $k\ge 1,$ define
\begin{align*}
  &S_{k+1}=\inf\{n:X_n>M_k\},\\
   &R_{k+1}=D\circ \theta_{S_{k+1}}+S_{k+1},\\
   & M_{k+1}=\sup\{X_n:0\le n\le R_{k+1}\}.
\end{align*}
We have $$0<S_1\le R_1\le S_2\le...\le S_{k}\le R_k \le \infty.$$
Let $$K=\inf\{k\ge 1: S_k<\infty \text{ but } R_k=\infty\}.$$ Let $\tau_0=0$ and whenever $K<\infty$ set $\tau_1=S_K.$
Define also $$\tau_2=[\tau_1](X.)+[\tau_1](X_{\tau_1+\cdot}-X_{\tau_1})$$ and recursively for $k\ge 2,$
$$\tau_{k+1}=[\tau_1](X.)+[\tau_k](X_{\tau_1+\cdot}-X_{\tau_1}).$$
Here and throughout, we use $[Y](Z.)$ to denote the random variable $Y$ defined by the process $Z.,$ and when $Y$ is defined by $X.,$ some times we write $[Y](X_.)$ simply as $Y.$
 \begin{lemma}\label{df}
   If $P(\lim_{n\rto}X_n=\infty)=1,$ then $P(D=\infty)>0.$
 \end{lemma} \proof Assume by contradiction that $P(D=\infty)=0.$ Then $P(D<\infty)=1.$ Consequently, $\mbp$-a.s., $P_\omega(D<\infty)=1.$ Therefore $P(\liminf_{n\rto}X_n\le 0)$ which  contradicts $P(\lim_{n\rto}X_n=\infty)=1.$\qed
 \begin{lemma}\label{kf}
   If $P(D=\infty)>0,$ then $P(K<\infty)=1$ and hence $P(\tau_1<\infty)=1.$
    \end{lemma}
    \proof For $k\ge1,$
    \begin{align*}
      P(R_k&<\infty)=\mbe\left(E_\omega\z(S_k<\infty, P^{X_{S_k}}_\omega(D<\infty)\y)\right)\\
       &=\sum_{x\in \mbz}\mbe\left(P_\omega(S_k<\infty, X_{S_k}=x)P^{x}_\omega(D<\infty)\right) \\
       &=\sum_{x\in \mbz}P(S_k<\infty, X_{S_k}=x)P(D<\infty) \\
          &=P(S_k<\infty)P(D<\infty) \\
         &\le P(R_{k-1}<\infty)P(D<\infty),
          \end{align*}
    where the third equality holds because $P_\omega(S_k<\infty, X_{S_k}=x)$ is $\sigma(\omega_i,i<x)$-measurable while $P^{x}_\omega(D<\infty)$ is $\sigma(\omega_i,i\ge x)$-measurable, and the environment is stationary under $\mbp.$
        By induction, we have that,  for $k\ge 1,$
    $$P(R_k<\infty)\le P(D<\infty)^k.$$
    Then an application of Borel-Cantelli lemma yields that $$P(R_k=\infty \text{ for some }k)=1.$$
    Consequently, $P(K<\infty)=1.$\qed

    Since $\gamma_2>0,$ $P(\lim_{n\rto}X_n=\infty)=1.$ Therefore we have from Lemma \ref{df} and Lemma \ref{kf} that $P(D=\infty)>0$ and hence $P(K<\infty)=P(\tau_1<\infty)=1.$
 For $k\ge 1,$ let $$\mathcal G_k=\sigma(\tau_1,...,\tau_k; (X_{\tau_k\wedge n})_{n\ge0},(\omega_{y})_{y< X_{\tau_k}}).$$ Consequently, the following lemma follows similarly as Kesten \cite{kes77} or Sznitman and Zerner \cite{sz99}.

    \begin{lemma}\label{tid} Suppose that condition (\ref{ell}) holds and $\gamma_2>0.$ Then for $k\ge 1,$ \begin{align*}&P\z[(X_{\tau_k+n}-X_{\tau_k})_{n\ge0}\in\cdot, (\omega_{y})_{y\ge X_{\tau_k}}\in \cdot|\mathcal G_k\y]\\
    &\quad\quad\quad\quad\quad\quad=P\z[(X_n)_{n\ge 0}\in\cdot, (\omega_{y})_{y\ge 0}\in\cdot|D=\infty\y].\end{align*}

\end{lemma}

 For $n\ge 0,$ define $T_n=\inf\{k:X_k> n\}.$

 \begin{lemma}\label{etau}
 Suppose that condition (\ref{ell}) holds and $\gamma_2>0.$ Then we have that \begin{equation}\label{taa}E(X_{\tau_1})=E(X_{S_1})+\frac{E(M+\xi; D<\infty)}{P(D=\infty)}<\infty,\end{equation} where $\xi:=X_{T_M}-M.$
    \end{lemma}
 \proof Note that
 \begin{align*}
   E(X_{\tau_1})&=\sum_{k\ge1} E(X_{\tau_1};K=k)=E(X_{S_k};S_k<\infty,D\circ \theta_{S_k}=\infty)\\
   &=\sum_{k\ge1}\sum_{x\in\mbz}\mbe\z[E_\omega(X_{S_k};S_k<\infty,X_{S_k}=x, D\circ \theta_{S_k}=\infty)\y]\\
 &=\sum_{k\ge1}\sum_{x\in\mbz}\mbe\z[xP^x_\omega( D=\infty)P_\omega(S_k<\infty, X_{S_k}=x)\y]
 \end{align*} where the third equality follows from the strong Markov property.
 Since $P^x_\omega( D=\infty)$ is $\sigma\{\omega_i:i\ge x\}$-measurable and $P_\omega(S_k<\infty, X_{S_k=x})$ is $\sigma\{\omega_i:i< x\}$-measurable, it follows by stationarity and independence of the environment that
 \begin{align*}
   E(X_{\tau_1})&=\sum_{k\ge1}\sum_{x\in\mbz}xP( D=\infty)P(S_k<\infty, X_{S_k=x})\\
   &=\sum_{k\ge1}P( D=\infty)E(X_{S_k};S_k<\infty).
  \end{align*}
 For $k\ge 2,$ we have that
 \begin{align*}
  E(X_{S_k};&S_k<\infty)=E(X_{S_k};S_{k-1}<\infty, D\circ \theta_{S_{k-1}}<\infty)\\
  &=\sum_{x\in\mbz}\mbe\z[E^x_\omega( x+M+\xi;D<\infty)P_\omega(S_{k-1}<\infty, X_{S_{k-1}}=x)\y]\\
  &\text{using again the stationarity and independence}\\
  &= \sum_{x\in\mbz}E( x+M+\xi;D<\infty)P(S_{k-1}<\infty, X_{S_{k-1}}=x)\\
&= E(X_{S_{k-1}};S_{k-1}<\infty)P(D<\infty)+E(M+\xi;D<\infty)P(S_{k-1}<\infty)
            \end{align*}
  Since $P(S_1<\infty)=1,$ by induction, it follows that
 \begin{align*}
   P(S_k<\infty)&=\sum_{x\in\mbz}\mbe\z[E^x_\omega(D<\infty)P_\omega(S_{k-1}<\infty,X_{S_{k-1}}=x)\y]\\
   &=P(S_{k-1}<\infty)P(D<\infty)=P(D<\infty)^{k-1}.
 \end{align*}
 Consequently, using again the induction, we have that
 \begin{align*}
  E(&X_{S_k};S_k<\infty)= E(X_{S_{k-1}};S_{k-1}<\infty)P(D<\infty)\\
&\quad\quad\quad\quad\quad+E(M+\xi;D<\infty)P(D<\infty)^{k-2}\\
&=E(X_{S_{k-2}};S_{k-2}<\infty)P(D<\infty)^2+2E(M+\xi;D<\infty)P(D<\infty)^{k-2}\\
&=E(X_{S_{1}})P(D<\infty)^{k-1}+(k-1)E(M+\xi;D<\infty)P(D<\infty)^{k-2}.
            \end{align*}
It follows from Theorem \ref{mt} that $E(M)<\infty.$
Therefore, by using the facts  $P(D=\infty)>0$ and $P(\xi=1)+P(\xi=2)=1,$ it could be concluded that
\begin{align*}
   E(X_{\tau_1})&=\sum_{k\ge1}P( D=\infty)E(X_{S_k};S_k<\infty)\\
   &=E(X_{S_1})+\frac{E(M+\xi; D<\infty)}{P(D=\infty)}<\infty.
  \end{align*}
The lemma is proved. \qed
\begin{remark}
  From the proof we see that the result of Lemma \ref{etau} could be strengthened. Indeed we have that $E(e^{c_3X_{\tau_1}})<\infty,$ for some $c_3>0.$
\end{remark}

Recall that $K=\inf\{k>0:S_k<\infty,D\circ\theta_{S_k}=\infty\}$ and $\tau_1=S_K.$
 Define $$\nu_1=\inf\{k\ge1:[X_{\tau_1}-M_{K-1}](X_{\tau_{k-1}+\cdot}-X_{\tau_{k-1}})=2\}.$$
 And for $k\ge 1.$ define recursively $$\nu_{k+1}=[\nu_1](X.)+[\nu_k](X_{\nu_1+\cdot}-X_{\nu_1}).$$

 \begin{lemma}\label{nid}
 Suppose that condition (\ref{ell}) holds and $\gamma_2>0.$ Then for $k\ge1,$
 \begin{align*}
   &P\z[\begin{array}{l}
          (X_{\tau_{\nu_k}+n}-X_{\tau_{\nu_k}})_{n\ge0}\in \cdot, (\omega_y)_{y\ge X_{\tau_{\nu_k}}}\in \cdot,\\
           (X_n)_{n\le \tau_{\nu_k}}\in\cdot, (\omega_y)_{y<X_{\tau_{\nu_k}}}\in\cdot
        \end{array}
\y]\\
   &\quad=P\z[(X_n)_{n\ge0}\in \cdot,(\omega_y)_{y\ge 0}\in \cdot|D=\infty \y] P\z[(X_n)_{n\le \tau_{\nu_k}}\in\cdot, (\omega_y)_{y<X_{\tau_{\nu_k}}}\in\cdot\y]
 \end{align*}
    \end{lemma}
 \proof Note that \begin{align*}
   &P\z[\begin{array}{l}
          (X_{\tau_{\nu_k}+n}-X_{\tau_{\nu_k}})_{n\ge0}\in \cdot, (\omega_y)_{y\ge X_{\tau_{\nu_k}}}\in \cdot,\\
           (X_n)_{n\le \tau_{\nu_k}}\in\cdot, (\omega_y)_{y<X_{\tau_{\nu_k}}}\in\cdot
        \end{array}
\y]\\
   &\quad\quad =\sum_{j=1}^\infty P\z[\begin{array}{l}
          (X_{\tau_{j}+n}-X_{\tau_{j}})_{n\ge0}\in \cdot, (\omega_y)_{y\ge X_{\tau_j}}\in \cdot,\\
           (X_n)_{n\le \tau_{j}}\in\cdot, (\omega_y)_{y<X_{\tau_{j}}}\in\cdot, \nu_k=j
        \end{array}\y].
           \end{align*}
   Since the event $\{\nu_k=j\}$ is $\mathcal G_j$-measurable, it follows from Lemma \ref{tid} that the right-hand side of the above equation equals to
   \begin{align*}
     &P\z[(X_n)_{n\ge0}\in \cdot,(\omega_y)_{y\ge 0}\in \cdot|D=\infty \y]  \\
     &\quad\quad \quad\quad\quad\times \sum_{j=1}^\infty P\z[(X_n)_{n\le \tau_{\nu_k}}\in\cdot, (\omega_y)_{y<X_{\tau_{\nu_k}}}\in\cdot,\nu_k=j\y]\\
     &\quad\quad=P\z[(X_n)_{n\ge0}\in \cdot,(\omega_y)_{y\ge 0}\in \cdot|D=\infty \y]  \\
     &\quad\quad \quad\quad\quad\times  P\z[(X_n)_{n\le \tau_{\nu_k}}\in\cdot, (\omega_y)_{y<X_{\tau_{\nu_k}}}\in\cdot\y]
   \end{align*}
   The lemma is proved. \qed

 Write $\nu=\nu_1$ for simplicity.

 \begin{lemma}\label{etnu}
 Suppose that condition (\ref{ell}) holds and $\gamma_2>0.$ Then \begin{equation}\label{ta}E(X_{\tau_{\nu}})=E(X_{\tau_1})+\frac{P(X_{\tau_1}-M_{K-1}=1)E(X_{\tau_1}| D=\infty)}{P(X_{\tau_1}-M_{K-1}=2|D=\infty)}<\infty.\end{equation}

 \end{lemma}
 \proof
Note that \begin{equation*}
   \{\nu=1\}=\{[X_{\tau_1}-M_{K-1}](X_.)=2\}
  \end{equation*} and for $l\ge2$
 \begin{equation*}
   \{\nu=l\}=\z\{\begin{array}{l}
                   [X_{\tau_1}-M_{K-1}](X_{\tau_{j}+\cdot}-X_{\tau_j})=1,\text{ for }j=0,...,l-2, \\
                  \text{ but }[X_{\tau_1}-M_{K-1}](X_{\tau_{l-1}+\cdot}-X_{\tau_{l-1}})=2
                 \end{array}\y\}.
  \end{equation*}
  For $j=0,...,l-1,$ denote simply $Y_j=[X_{\tau_1}-M_{K-1}](X_{\tau_{j}+\cdot}-X_{\tau_j}).$ It is easy to see that $Y_j$ is $\mathcal G_{j+1}$-measurable and it follows from Lemma \ref{tid} that $Y_0, Y_1,...,Y_{l-1}$ share the same distribution $$P\z([X_{\tau_1}-M_{K-1}](X.)\in \cdot|D=\infty\y)$$ except for $Y_0.$ In the proof of the lemma, we write temporarily $$p=P\z(Y_0=2|D=\infty\y)$$ for simplicity. By total probability
   we have that \begin{equation}\label{tnu}\begin{split}
   E(X_{\tau_{\nu}})&=\sum_{l=1}^\infty E\z(X_{\tau_{l}};\nu=l\y)=\sum_{l=1}^\infty E\z(\sum_{i=1}^l X_{\tau_{i}}-X_{\tau_{i-1}};\nu=l\y)\\
   &=\sum_{l=1}^\infty E\z(\sum_{i=1}^l (X_{\tau_{i}}-X_{\tau_{i-1}})\mathbf{1}_{Y_{l-1}=2}\prod_{j=0}^{l-2}\mathbf{1}_{Y_j=1}\y)
 \end{split}\end{equation}
 with the convention that the empty product equals to $1.$

 On one hand, fixing $l\ge3,$ since $\prod_{j=0}^{l-2}\mathbf{1}_{Y_j=1}$ is $\mathcal G_{l-1}$-measurable, then by Lemma \ref{tid}
 \begin{align*}
   &E\z((X_{\tau_{l}}-X_{\tau_{l-1}})\mathbf{1}_{Y_{l-1}=2}\prod_{j=0}^{l-2}\mathbf{1}_{Y_j=1}\y)\\
   &\quad\quad =E\z(E\z[(X_{\tau_{l}}-X_{\tau_{l-1}})\mathbf{1}_{Y_{l-1}=2}\prod_{j=0}^{l-2}\mathbf{1}_{Y_j=1}\Big |G_{l-1}\y]\y)\\
   &\quad\quad =E(X_{\tau_{1}}\mathbf{1}_{Y_{0}=2}|D=\infty) E\z[\prod_{j=0}^{l-2}\mathbf{1}_{Y_j=1}\y].
   \end{align*}
    Using the same trick for  $l-1$  times, we have that
    \begin{equation}\label{sa}\begin{split}
   &E\z((X_{\tau_{l}}-X_{\tau_{l-1}})\mathbf{1}_{Y_{l-1}=2}\prod_{j=0}^{l-2}\mathbf{1}_{Y_j=1}\y)\\
   &\quad\quad =E(X_{\tau_{1}}\mathbf{1}_{Y_{0}=2}|D=\infty) P(Y_0=1)(1-p)^{l-2}.
\end{split}
\end{equation}
 Similarly as above, for $2\le i\le l-1,$ we have that
 \begin{equation}\label{sb}\begin{split}
   &E\z((X_{\tau_{i}}-X_{\tau_{i-1}})\mathbf{1}_{Y_{l-1}=2}\prod_{j=0}^{l-2}\mathbf{1}_{Y_j=1}\y)\\
   &\quad\quad =E\z(X_{\tau_{1}}\mathbf{1}_{Y_{0}=1}|D=\infty\y)P(Y_0=1)(1-p)^{l-3}p
 \end{split}
\end{equation}
  and   \begin{equation}\label{sc}\begin{split}
   E\z((X_{\tau_{1}}-X_{\tau_{0}})\mathbf{1}_{Y_{l-1}=2}\prod_{j=0}^{l-2}\mathbf{1}_{Y_j=1}\y) =E\z(X_{\tau_{1}}\mathbf{1}_{Y_{0}=1}\y)(1-p)^{l-2}p.
\end{split}
\end{equation}
 On the other hand, for $l=1$
 \begin{equation}\label{sd}\begin{split}
   E\z(\sum_{i=1}^l (X_{\tau_{i}}-X_{\tau_{i-1}})\mathbf{1}_{Y_{l-1}=2}\prod_{j=0}^{l-2}\mathbf{1}_{Y_j=1}\y)=E(X_{\tau_1}\mathbf{1}_{Y_0=2})
 \end{split}\end{equation}
  while for $l=2$
    \begin{equation}\label{se}\begin{split}
   &E\z(\sum_{i=1}^l (X_{\tau_{i}}-X_{\tau_{i-1}})\mathbf{1}_{Y_{l-1}=2}\prod_{j=0}^{l-2}\mathbf{1}_{Y_j=1}\y)\\
   &\quad\quad=p E(X_{\tau_1}\mathbf1_{Y_0=1})+P(Y_0=1)E(X_{\tau_1}\mathbf1_{Y_0=2}|D=\infty).
 \end{split}\end{equation}
 Substituting (\ref{sa})-(\ref{se}) to (\ref{tnu}), we have that
   \begin{equation*}\begin{split}
   E(X_{\tau_{\nu}})&=E(X_{\tau_1}\mathbf{1}_{Y_0=2})+p E(X_{\tau_1}\mathbf1_{Y_0=1})+P(Y_0=1)E(X_{\tau_1}\mathbf1_{Y_0=2}|D=\infty)\\
 &\quad\quad+\sum_{l=3}^\infty (l-2)E\z(X_{\tau_{1}}\mathbf{1}_{Y_{0}=1}|D=\infty\y)P(Y_0=1)(1-p)^{l-3}p \\
  &\quad\quad\quad\quad\quad\quad+E\z(X_{\tau_{1}}\mathbf{1}_{Y_{0}=1}\y)(1-p)^{l-2}p\\
   &\quad\quad\quad\quad\quad\quad+E\z(X_{\tau_{1}}\mathbf{1}_{Y_{0}=2}|D=\infty\y)P(Y_0=1)(1-p)^{l-2}.
 \end{split}\end{equation*}
  Consequently,
  $$ E(X_{\tau_{\nu}})=E(X_{\tau_1})+\frac{P(X_{\tau_1}-M_{K-1}=1)E(X_{\tau_1}|D=\infty)}{P(X_{\tau_1}-M_{K-1}=2|D=\infty)}.$$
  To show the finiteness of $E(X_{\tau_{\nu}}),$  note that \begin{equation*}\label{xmg}\begin{split}
    P(X_{\tau_1}&-M_{K-1}=2,D=\infty)\ge P(X_{S_1}-M_0=2, K=1, D=\infty)\\
    &=\mbe(\omega_0(2)P_\omega^2(D=\infty))=\mbe(\omega_0(2))P(D=\infty)>0.
  \end{split}\end{equation*}
 But by Lemma \ref{etau}, $E(X_{\tau_1})$ is finite. Then the lemma follows.\qed

We are know ready to prove Theorem \ref{main}.

\noindent{\bf Proof of Theorem \ref{main}.} Recall that $\tau_k,k\ge 1$ are successive renewal epochs and especially, each of $\nu_k,k\ge 1$  is a renewal epoch such that $X_{\nu_k}$ minuses the maximum before $\tau_{\nu_k}$ equal to $2.$ Note that after $\tau_{\nu_k}$ the walk will not visit any site in the left side of $X_{\tau_{\nu_k}}$ any longer. Therefore the walk never visits the site $X_{\tau_{\nu_k}}-1.$  And all  other sites in $[X_{\tau_{\nu_{k-1}}},X_{\tau_{\nu_k}-1}]$ were visited at least once, since the downward jumps of the walk are nearest neighbor. We conclude that there are rightly $k$ sites, say $X_{\tau_{\nu_i}}-1,\ i=1,..,k$ in $[0,X_{\tau_{\nu_k}}]$ which are never visited by the walk.

Let $$N(x)=\#\{S_n: 0\le S_n\le x, n\ge 0\}.$$
For $x\in \mbz_+,$ there exists a unique random number $k(x)$ such that $$X_{\tau_{\nu_k}}\le x<X_{\tau_{\nu_{k+1}}}.$$
We have that \begin{align}\label{jba}
  \frac{N(X_{\tau_{\nu_k}})}{X_{\tau_{\nu_{k+1}}}}\le \frac{N(x)}{x}\le\frac{N(X_{\tau_{\nu_{k+1}}})}{X_{\tau_{\nu_{k}}}}.
\end{align}
But one follows from Lemma \ref{nid} that $$(X_{\tau_{\nu_1}},\tau_{\nu_1}),(X_{\tau_{\nu_2}}-X_{\tau_{\nu_{1}}},\tau_{\nu_2}-\tau_{\nu_{1}}),...,(X_{\tau_{\nu_k}}-X_{\tau_{\nu_{k-1}}},\tau_{\nu_k}-\tau_{\nu_{k-1}}),...$$
are independent and $$(X_{\tau_{\nu_2}}-X_{\tau_{\nu_{1}}},\tau_{\nu_2}-\tau_{\nu_{1}}),...,(X_{\tau_{\nu_k}}-X_{\tau_{\nu_{k-1}}},\tau_{\nu_k}-\tau_{\nu_{k-1}}),...$$
are all distributed as $(X_{\tau_{\nu_1}},\tau_{\nu_1})$ under $P(\cdot|D=\infty).$ Thus an application of the strong law of large numbers yields that $P$-a.s.,
\begin{align}\label{jbb}
 \lim_{k\rto}\frac{N(X_{\tau_{\nu_{k}}})}{X_{\tau_{\nu_{k+1}}}}=\lim_{k\rto}\frac{X_{\tau_{\nu_{k}}}-k}{X_{\tau_{\nu_{k+1}}}}=\frac{E(X_{\tau_{\nu}}|D=\infty)-1}{E(X_{\tau_{\nu}}|D=\infty)}.
\end{align}
Similarly we have that $P$-a.s., \begin{align}\label{jbc}
 \lim_{k\rto}\frac{N(X_{\tau_{\nu_{k+1}}})}{X_{\tau_{\nu_{k}}}}=\frac{E(X_{\tau_{\nu}}|D=\infty)-1}{E(X_{\tau_{\nu}}|D=\infty)}.
\end{align}
Note that by Lemma \ref{etnu}, $E(X_{\tau_{\nu}}|D=\infty)<\infty.$ Then it follows from (\ref{jba}), (\ref{jbb}) and (\ref{jbc}) that $P$-a.s.,\begin{equation}\label{tb}
  \lim_{n\rto}\frac{N(x)}{x}=\frac{E(X_{\tau_{\nu}}|D=\infty)-1}{E(X_{\tau_{\nu}}|D=\infty)}<1.
\end{equation}
Theorem \ref{main} is proved.\qed

\section{Proof of Theorem \ref{mt}}
For $n\ge 0,$ define $\tilde T_n=\inf\{k\ge0: X_k=n\}.$  We have that
\begin{align*}
 P(M=n,& D<\infty)=\mbe(P_\omega(M=n,D<\infty))\\
 &=\mbe\z(P_\omega(M=n,\tilde T_n<\infty, D<\infty)\y)\\
 &\le \mbe\z(P^n_\omega(\text{the walk hits }(-\infty,0)\text{ before }[n+1,\infty)\y)\\
&= \mbe\z(P_\omega(\text{the walk hits }(-\infty,-n)\text{ before }[1,\infty)\y).
\end{align*}
By some delicate calculation, \begin{align*}
  P_\omega[&\text{the walk hits }(-\infty,-n)\text{ before }[1,\infty)]\\
  &=\frac{1}{1+\sum_{j=-n}^0e_1A_jA_{j+1}\cdots A_0e_1^t}
\end{align*}
referring to \cite{brb} for the details.
Then it follows that
\begin{align*}
 P(&M=n,D<\infty)=\mbe(P_\omega(M=n,D<\infty))\\
 &\le\mbe\z( \frac{1}{1+\sum_{j=-n+1}^0e_1A_jA_{j+1}\cdots A_0e_1^t}\y)\\
  &\le \mbe\z(\frac{1}{1+e_1A_{-n+1}\cdots A_0e_1^t}\y)\\
  &\le \mbe\z(\frac{1}{1+\left(a_{-n+1}(1)e_1+a_{-n+1}(2)e_2\right)A_{-n+2}\cdots A_{-1}\z(a_0(1)e_1+e_2\y)^t}\y).
 \end{align*}
 Using the stationarity of the environment, we have that
 \begin{equation}\label{anor}\begin{split}
 P(&M=n,D<\infty)\\
 &\le \mbe\z(\frac{1}{1+\min\{a_{0}(1),a_{0}(2)\}\min\{a_{n-1}(1),1\}\mathbf1A_{1}\cdots A_{n-2}\mathbf1^t}\y)\\
 &\le \mbe\z(\frac{1}{1+\min\{a_{0}(1),a_{0}(2)\}\min\{a_{n-1}(1),1\}\|A_{1}\cdots A_{n-2}\|}\y)
 \end{split}\end{equation}
 where $\mathbf 1=e_1+e_2.$

Under condition (\ref{ell}), $A_i, i\in\mbz$ are all positive matrices satisfying $$A_iA_{i+1}\gg 0,$$ where a matrix $A\gg0$ means all its entries  are strictly positive. Then it follows  from Frobenious theory that for each $A_i$ there is a single eigenvalue $\lambda$ which dominates all other eigenvalues. Therefore, the support of $\mbp$ is contracting. Also by the special structures of $M_i,\ i\in \mathbb Z,$  the group generated by the support of $\mbp$ is strong irreducible. Hence a large deviation argument (see \cite{bl85}) yields that for some   small $\delta>0$ there exists for  $\eta\in(\gamma_2-\delta,\gamma_2)$  a number $c_\eta>0$ such that
$$\lim_{n\rto}\frac{1}{n}\log\mbp\z(\frac{1}{n}\log\|A_1\cdots A_{n-2}\|<\eta\y)=-c_\eta.$$
 Write temporarily $B_0(\omega):=\min\{a_{0}(1),a_{0}(2)\}$ and $B_{n-1}(\omega):= \min\{a_{n-1}(1),1\}$ for simplicity. By (\ref{ell})
 $\mbe(|\log |B_0(\omega)|+|\log B_{n-1}(\omega)|)<\infty.$
Hence we have that \begin{align*}
  \lim_{n\rto}\frac{1}{n}\log\mbp\z[\frac{1}{n}\log\z\{ B_0(\omega) B_{n-1}(\omega)\|A_1\cdots A_{n-2}\|\y\}<\rho\y]=-c_4.
\end{align*}
  for some $0<\rho<\gamma_2$ and $c_4>0.$
Fix $0<c_5<c_4,$ there exist $N>0$ such that for all $n>N,$
 $$ \mbp\z[\frac{1}{n}\log\z\{ B_0(\omega) B_{n-1}(\omega)\|A_1\cdots A_{n-2}\|\y\}<\rho\y]<e^{- c_5n}.$$
Substituting the above inequality to (\ref{anor}), we have that
\begin{equation*}\begin{split}
 P(&M=n,D<\infty)\\
 &\le e^{-\rho n} +\mbp\z( B_0(\omega) B_{n-1}(\omega)\|A_1\cdots A_{n-2}\|<e^{\rho n}\y)\\
 &\le c_6e^{-c_7n}
 \end{split}\end{equation*}
for some $c_6,c_7>0.$
Consequently, for some $c_1,c_2>0,$
$$ P(M>n,D<\infty)<c_1e^{-c_2n}.$$ Theorem \ref{mt} is proved.\qed

\noindent{\large{\bf \large Acknowledgements:}} The author was grateful to Prof. Vladimir Vatutin for his help when writing the paper.


\end{document}